\pdfoutput=1
\RequirePackage{ifpdf}
\ifpdf 
\documentclass[pdftex]{sigma}
\else
\documentclass{sigma}
\fi

\newtheorem{prob}{Problem}

{\theoremstyle{definition}
\newtheorem{re}{Remark}}

\usepackage{esint}

\DeclareMathOperator{\Area}{Area}

\DeclareMathOperator{\Haw}{H}

\def\Hm{\mathfrak{m}_{\scriptscriptstyle\Haw}}

\begin{document}
\allowdisplaybreaks

\newcommand{\arXivNumber}{2001.05607}

\renewcommand{\thefootnote}{}

\renewcommand{\PaperNumber}{030}

\FirstPageHeading

\ShortArticleName{NNSC-Cobordism of Bartnik Data in High Dimensions}

\ArticleName{NNSC-Cobordism of Bartnik Data\\ in High Dimensions\footnote{This paper is a~contribution to the Special Issue on Scalar and Ricci Curvature in honor of Misha Gromov on his 75th Birthday. The full collection is available at \href{https://www.emis.de/journals/SIGMA/Gromov.html}{https://www.emis.de/journals/SIGMA/Gromov.html}}}

\Author{Xue HU~$^\dag$ and Yuguang SHI~$^\ddag$}

\AuthorNameForHeading{X.~Hu and Y.~Shi}

\Address{$^\dag$~Department of Mathematics, College of Information Science and Technology,\\
\hphantom{$^\dag$}~Jinan University, Guangzhou, 510632, P.R.~China}
\EmailD{\href{mailto:thuxue@jnu.edu.cn}{thuxue@jnu.edu.cn}}

\Address{$^\ddag$~Key Laboratory of Pure and Applied Mathematics, School of Mathematical Sciences,\\
\hphantom{$^\ddag$}~Peking University, Beijing, 100871, P.R.~China}
\EmailD{\href{mailto:ygshi@math.pku.edu.cn}{ygshi@math.pku.edu.cn}}

\ArticleDates{Received January 22, 2020, in final form April 13, 2020; Published online April 20, 2020}

\Abstract{In this short note, we formulate three problems relating to nonnegative scalar curvature (NNSC) fill-ins. Loosely speaking, the first two problems focus on: {\it When are $(n-1)$-dimensional Bartnik data $\big(\Sigma_i ^{n-1}, \gamma_i, H_i\big)$, $i=1,2$, NNSC-cobordant}? (i.e., there is an $n$-dimensional compact Riemannian manifold $\big(\Omega^n, g\big)$ with scalar curvature $R(g)\geq 0$ and the boundary $\partial \Omega=\Sigma_{1} \cup \Sigma_{2}$ such that $\gamma_i$ is the metric on $\Sigma_i ^{n-1}$ induced by~$g$, and~$H_i$ is the mean curvature of $\Sigma_i$ in $\big(\Omega^n, g\big)$). If $\big(\mathbb{S}^{n-1},\gamma_{\rm std},0\big)$ is positive scalar curvature (PSC) cobordant to  $\big(\Sigma_1 ^{n-1}, \gamma_1, H_1\big)$, where $\big(\mathbb{S}^{n-1}, \gamma_{\rm std}\big)$ denotes the standard round unit sphere then $\big(\Sigma_1 ^{n-1}, \gamma_1, H_1\big)$ admits an NNSC fill-in. Just as Gromov's conjecture is connected with positive mass theorem, our problems are connected with Penrose inequality, at least in the case of $n=3$. Our third problem is on $\Lambda\big(\Sigma^{n-1}, \gamma\big)$ defined below.}


\Classification{53C20; 83C99}

\begin{flushright}
\begin{minipage}{80mm}
\it Dedicate this paper to Professor Misha Gromov \\
on the occasion of his 75th birthday.
\end{minipage}
\end{flushright}

\renewcommand{\thefootnote}{\arabic{footnote}}
\setcounter{footnote}{0}

Bartnik data $\big(\Sigma^{n-1}, \gamma, H\big)$ consists of an $(n-1)$-dimensional orientable Riemannian mani\-fold~$\big(\Sigma^{n-1}, \gamma\big)$ and a smooth function $H$ defined on~$\Sigma^{n-1}$ which serves as the mean curvature of~$\Sigma^{n-1}$. One basic problem in Riemannian geometry is to study: {\it under what conditions is it that~$\gamma$ is induced by a Riemannian metric $g$ with nonnegative scalar curvature, for example, defined on $\Omega^n$, and $H$ is the mean curvature of $\Sigma$ in $\big(\Omega^n, g\big)$ with respect to the outward unit normal vector?} Indeed, this problem was proposed by M.~Gromov recently (see \cite[Problem~A]{Gromov2} and \cite[Sections~3.3 and~3.6]{Gromov4}).

 On the other hand, when $n=3$, for each Bartnik data $\big(\Sigma^2, \gamma, H\big)$ may be associated with certain quasi-local masses, for instance, when the Gaussian curvature~$K$ of $\gamma$ is positive, $(\mathbb{S}^2, \gamma)$ can be isometrically embedded into $\mathbb{R}^3$ with mean curvature $H_0$ (with respect to the outward unit normal vector of the embedded image in $\mathbb{R}^3$), with this embedding we may define Brown--York mass for $\big(\mathbb{S}^2, \gamma, H\big)$ \cite{BrownYork-1, BrownYork-2} as
\begin{gather*}
 \mathfrak{m}_{\rm BY}\big(\mathbb{S}^2; \gamma, H\big)=\frac1{8\pi}\int_{\mathbb{S}^2}(H_0-H) \, {\rm d}\sigma_\gamma.
\end{gather*}

 If $\big(\mathbb{S}^2, \gamma, H\big)$ admits an NNSC fill-in and $H>0$, it was shown that $\mathfrak{m}_{\rm BY}\big(\mathbb{S}^2; \gamma, H\big)\geq 0$~\cite{ST1}. There are several pieces of interesting work on NNSC fill-ins relating to positivity of Brown--York mass (for instance see~\cite{J,JMT}). Obviously, positivity of Brown--York mass is one necessary condition for the existence of such a fill-in, but it is far from sufficient. {\it It was shown that for Bartnik data $\big(\mathbb{S}^2, \gamma, H\big)$ with positive Gaussian curvature and $H>0$, let~$H_0$ be the mean curvature of isometric embedding of $\big(\mathbb{S}^2, \gamma\big)$ in~$\mathbb{R}^3$, if $\mathfrak{m}_{\rm BY}\big(\mathbb{S}^2; \gamma, H\big)=0$ and $H\neq H_0$ then there is a constant $\epsilon$ depending only on $\big(\mathbb{S}^2, \gamma, H\big)$ such that for any $\tilde {H}>H-\epsilon$, $\big(\mathbb{S}^2, \gamma, \tilde {H}\big)$ admits no NNSC fill-ins} \cite[Theorem~3]{JMT}.

 If $K>-\kappa^2$ where $\kappa$ is a constant, then $\big(\mathbb{S}^2, \gamma\big)$ can be isometrically embedded into the hyperbolic space with constant sectional curvature $-\kappa^2$, and we can make use of such embedding to define a generalized Brown--York mass, moreover if $H>0$ we were able to prove its positivity~\cite{ShiTam2007}. Clearly, this positivity of generalized Brown--York mass is also a kind of necessary condition for the Bartnik data with $K>-\kappa^2$ and $H>0$ to admit NNSC fill-ins.

 For Bartnik data $\big(\Sigma^2, \gamma, H\big)$, we can define its Hawking mass as following:
\begin{gather*}\Hm(\Sigma, \gamma, H)=\sqrt{\frac{\Area(\Sigma)}{16\pi}}\left(1-\frac{1}{16\pi}\int_{\Sigma}H^2\,{\rm d}\sigma_\gamma\right).
\end{gather*}
 It should be interesting to explore similar relation between Hawking mass or other quasi-local masses of the Bartnik data with its NNSC fill-ins. Unfortunately, it is not easy to obtain a lower bound of the Hawking mass which depends only on $\big(\Sigma^2, \gamma\big)$.

 In the investigation of above Gromov's NNSC fill-in problem, we often need to deal with NNSC-cobordisms of Bartnik data which may have its own interests. More specifically, {\it given Bartnik data $\big(\Sigma^{n-1}_i, \gamma_i, H_i\big)$, $i=1,2$, we say $\big(\Sigma^{n-1}_1, \gamma_1, H_1\big)$ is NNSC-cobordant to $\big(\Sigma^{n-1}_2, \gamma_2, H_2\big)$ if there is an orientable $n$-dimensional manifold $\big(\Omega^n,g\big)$ with $\partial \Omega^n=\Sigma_1^{n-1} \cup \Sigma^{n-1}_2$, $R(g)\geq 0$, $\gamma_i=g|_{\Sigma_i}$, $i=1,2$, $H_1$ is the mean curvature of $\Sigma^{n-1}_1$ in $\big(\Omega^n,g\big)$ with respect to inward unit normal vector, and $H_2$ is the mean curvature of $\Sigma^{n-1}_2$ in $\big(\Omega^n,g\big)$ with respect to outward unit normal vector.} Our first problem is:
\begin{prob}\label{1}
	Given Bartnik data $\big(\Sigma^{n-1}_i, \gamma_i, H_i\big)$, $i=1,2$, when are they NNSC-cobordant?
\end{prob}

By using surgery arguments (see \cite{GL,SY1}), it is not difficult to show that if Bartnik data $\big(\Sigma^{n-1}_i\!, \gamma_i, H_i\big)$, $i=1,2$ can be filled in with positive scalar curvature metrics, then $\big(\Sigma^{n-1}_1\!, \gamma_1, -H_1\big)\!$ is NNSC-cobordant to $\big(\Sigma^{n-1}_2, \gamma_2, H_2\big)$. Another possible relevant notion to this is so called ``PSC-concordant''.
 Namely, {\it two PSC-metrics $\gamma_0$ and $\gamma_1$ on $\Sigma^{n-1} $ are said to be PSC-concordant if there is a PSC-metric $g$ on the cylinder $\Sigma \times I$ which are the product $\gamma_0 +{\rm d} t^2$ near $\Sigma \times \{0\}$ and $\gamma_1 +{\rm d}t^2$
 near $\Sigma\times \{1\}$} (see~\cite{Wa}), in that case, $\big(\Sigma^{n-1}, \gamma_0,0\big)$ is NNSC-cobordant to $\big(\Sigma^{n-1}, \gamma_1,0\big)$. By index theory, it is known that there are countable infinity distinct PSC-concordant classes for~$\mathbb{S}^{4k-1}$, for any positive integer $k\geq 2$. When two PSC-metrics $\gamma_0$ and $\gamma_1$ are isotopic, i.e., they can be connected by a continuous path $\gamma_t$, $t\in [0,1]$, and for each $t \in [0,1]$, $\gamma_t$ is a PSC-metric. Then we may use quasi-spherical metric to show that {\it if $H_1$ is not too large then $\big(\mathbb{S}^2, \gamma_0, H_0\big)$ is NNSC-cobordant to $\big(\mathbb{S}^2,\gamma_1, H_1\big)$, here $H_0$ can be any given smooth positive function} (see \cite{Bartnik, ST1,ST2}). On the other hand, {\it when $H_1$ is large enough we are able to show $\big(\mathbb{S}^2, \gamma_i, H_i\big)$, $i=0,1$, cannot be NNSC-cobordant}~\cite{BS1}.

 Let $\gamma_0$ be a Riemannian metric on~$\mathbb{S}^2$ with its first eigenvalue $\lambda_{1}(-\Delta_0+K)>0$, here $\Delta_0$ is the Laplacian operator of $\gamma_0$,
 then it was shown in~\cite{MS} {\it that $\big(\mathbb{S}^2,\gamma_0, 0 \big)$ is NNSC-cobordant to $\big(\mathbb{S}^2,\gamma_{\rm rou}, H \big)$ provided $\Hm\big(\mathbb{S}^2,\gamma_{\rm rou}, H\big)> \sqrt{\frac{{\rm Area}(\mathbb{S}^2,\gamma_0)}{16\pi}}$, here $\gamma_{\rm rou}$ denotes the round metric on~$\mathbb{S}^2$}. For a generalization to the case of Bartnik data with constant mean curvature surfaces see \cite[Theorem~1.1]{PCMM}, and higher-dimensional analogues see \cite[Theorems~1.1 and~1.2]{CPM}, and \cite[Proposition~2.1]{MM1}. An NNSC fill-in by a conformal blow-down argument which may have deep relation to Problem~\ref{1} please see the proof of Theorem~1.2 in~\cite{HM}. For deep discussion on PSC-concordant relation for two PSC-metrics on a manifold from topological point of view, please see~\cite{Wa2, Wa3} and references therein.

 As we mentioned above, one obstruction of the above NNSC fill-in problem is from positivity of certain quasi-local mass (for instance, Brown--York mass, see \cite{ST1, SWWZ}). It may be reasonable to think that there may be a potential obstruction of NNSC-cobordism problem which is from Penrose-type inequality (for Penrose inequality, see \cite{Bray1,HI}, for local Penrose inequality, see \cite{ LM,MX, ST3,SWY}). For instance, we observed that {\it if $\big(\mathbb{S}^2, \gamma_2, H_2\big)$ is with positive Gaussian curvature and $H_2>0$, and $\big(\Sigma_1^2, \gamma_1, H_1\big)$ is NNSC-cobordant to $\big(\mathbb{S}^2, \gamma_2, H_2\big)$, then $\mathfrak{m}_{\rm BY}\big(\mathbb{S}^2; \gamma_2, H_2\big)\geq \Hm\big(\Sigma_1^2, \gamma_1, H_1\big)$ provided $\Hm\big(\Sigma_1^2, \gamma_1, H_1\big)\leq 0$}~\cite{BS1}.

 To our knowledge, even the following simple case is still unknown:
\begin{prob}\label{2}Given Bartnik data $\big(\mathbb{S} ^{n-1}, g_1, H\big)$ and $\big(\mathbb{S} ^{n-1}, g_0, 0\big)$, both are with positive scalar curvature, what is the largest $\inf\limits_{\mathbb{S} ^{n-1}} H$ so that $\big(\mathbb{S} ^{n-1}, g_0, 0\big)$ is NNSC-cobordant to $\big(\mathbb{S} ^{n-1}, g_1, H\big)$?
\end{prob}
\begin{re}\quad
\begin{itemize}\itemsep=0pt
	\item By the arguments of \cite[Theorem~1.4]{SWWZ} and some gluing technique, we are able to show that for any PSC-metric $g_1$ on $\mathbb{S} ^{n-1}$, no matter whether~$g_1$ is PSC-concordant to~$g_0$ or not, there is a constant $H$ so that $\big(\mathbb{S} ^{n-1}, g_0, 0\big)$ is NNSC-cobordant to $\big(\mathbb{S} ^{n-1}, g_1, H\big)$ and the ambient manifold bounded by these Bartnik data is diffeomorphic to $\mathbb{S} ^{n-1}\times [0,1]$ provided~$g_0$ is the standard round metric on $\mathbb{S} ^{n-1}$ \cite{BS1}.
	\item If $g_0$ is the standard round metric on $\mathbb{S} ^{n-1}$, then by gluing arguments, the largest $\inf\limits_{\mathbb{S} ^{n-1}} H$ in Problem~\ref{2} is the corresponding number for $\big(\mathbb{S} ^{n-1}, g_1, H\big)$ to admit NNSC fill-ins.\footnote{We are grateful to the referee for pointing this fact to us.}
	
\item As we know, $\int_{\Sigma ^{2}}H \,{\rm d}\mu_1$ and $\int_{\Sigma ^{2}}H^2 \,{\rm d}\mu_1$ are closely related to Brown--York mass and Haw\-king mass respectively, they are also involved in classical Minkowski's inequality for a~convex surface and Willmore functional for a surface in $\mathbb{R}^3$, so, it may also be interesting to ask what the possible largest values of $\int_{\mathbb{S} ^{n-1}}H \,{\rm d}\mu_1$ and $\int_{\mathbb{S} ^{n-1}}H^2 \,{\rm d}\mu_1$ are, especially for $n=3$. 	
 \end{itemize}
\end{re}

For an orientable closed null-cobordant Riemannian manifold $\big(\Sigma^{n-1},\gamma\big)$, define $\Lambda\big(\Sigma^{n-1},\gamma\big)$ by
\begin{gather*}
\Lambda\big(\Sigma^{n-1},\gamma\big)=\sup\left\{\int_{\Sigma}H\,\mathrm d\mu_\gamma\,\Big |\ \big(\Sigma^{n-1},\gamma,H\big)\text{ admits an NNSC fill-in}\right\}.
\end{gather*}

In the case of $n=3$ and $H>0$, the above $\Lambda$ was introduced in~\cite{MM, MMT}, and also some interesting properties were discussed therein. An open problem on an estimate of $\Lambda\big(\Sigma^{n-1},\gamma\big)$ was proposed in \cite[p.~31]{Gromov4}, and a partial result in the case of $H>0$ was obtained in \cite[Theorem~1.3]{SWWZ}.

Suppose $\big(\mathbb{S}^2,\gamma\big)$ is a $2$-dimensional surface with positive Gaussian curvature, then it can be isometrically embedded into $\mathbb{R}^3$, let~$H_0$ be the mean curvature of the embedding image with respect to the outward unit normal vector, then we have:

\begin{prob}Is 	$ \Lambda\big(\mathbb{S}^2,\gamma\big)=\int_{\mathbb{S}^2}H_0 \,{\rm d} \mu_\gamma $?
\end{prob}

The affirmative answer implies the positivity of Brown--York mass without assumption of positivity of the mean curvature.

\subsection*{Acknowledgements}
The authors would like to thank Dr.~Georg Frenck for his comments to clarify some notions in this note, and are also deeply grateful to anonymous referees for invaluable suggestions on the exposition. The research of the first and the second author was partially supported by NSFC 11701215, NSFC 11671015 and 11731001 respectively.

\pdfbookmark[1]{References}{ref}
\LastPageEnding

\end{document}